\documentclass[%
 reprint,
 amsmath,amssymb,
 aps,
]{revtex4-2}

\usepackage{braket}
\usepackage{orcidlink}
\usepackage{graphicx}% Include figure files
\usepackage{dcolumn}% Align table columns on decimal point
\usepackage{bm}% bold math
\usepackage{comment}
\usepackage{amsmath}
\usepackage{caption}
\usepackage{subcaption}
\newcommand{\figsize}{0.23}

\begin{document}

\preprint{APS/123-QED}

\title{A practical recipe for variable-step finite differences via equidistribution}% Force line breaks with \\
%\thanks{A footnote to the article title}%

\author{Mário B. Amaro \orcidlink{0000-0003-3408-8533}}
 %\email{marioam@kth.se}
\affiliation{Department of Physics, Stockholm University, Stockholm S-10691, Sweden}%

\date{\today}

\begin{abstract}
We describe a short, reproducible workflow for applying finite differences on nonuniform grids determined by a positive weight function g. The grid is obtained by equidistribution, mapping uniform computational coordinates $\xi\in[0,1]$ to physical space by the cumulative integral $S(x)=\int_a^x\!1/g(s)\,ds$ and its inverse, and in multiple dimensions by the corresponding variable-diffusion (harmonic) mapping with tensor $P=(1/g)I$. We then use the standard three-point central stencils on uneven spacing for first and second derivatives. We collect the formulas, state the mild constraints on g (positivity, boundedness, integrability), and provide a small reference implementation. Finally, we solve the 2D time-independent Schrödinger equation for a harmonic oscillator on uniform vs. variable meshes, showing the expected improvement in resolving localized eigenfunctions without increasing matrix size. We intend this note as a how-to reference rather than a new method, consolidating widely used ideas into a single, ready-to-use recipe, claiming no novelty.
\end{abstract}

%Use showkeys class option if keyword
                              %display desired
\maketitle

%\tableofcontents

\section{Introduction}

The Finite Difference Method (FDM) is a family of techniques that is a cornerstone of numerical analysis, often used to approximate the solutions to partial differential equations (PDEs) in fields as different as electromagnetism \cite{kunz1993} \cite{Zhou1993}, fluid mechanics \cite{godunov1959} \cite{narasimhan1976} and financial mathematics \cite{duffy2013} \cite{HullWhite1990}. This method relies on discretizing the parameter space in a mesh and approximate derivatives using finite differences, which effectively reduces the problem to a linear algebra one, making it significantly simpler \cite{Grossmann2007}. With a sufficiently small step between gridpoints, we numerically approach the exact solution. A standard example is the one-dimensional uniform grid case. Here, we approximate the first and second order derivatives via the most simple finite differences
\begin{equation}
    f'_h(x) \approx \frac{f(x+h) - f(x-h)}{2h},
    \label{finitediff_order1}
\end{equation}
\begin{equation}
    f''_h(x) \approx \frac{f(x+h) - 2 f(x) + f(x-h)}{h^2},
    \label{finitediff_order2}
\end{equation}
where it is easily seen that, in the limit where h tends to zero, Equation \ref{finitediff_order1} corresponds to the definition of derivative in Calculus. For variable meshes, there are various techniques to extend this formulation for that case, and they are in general application-specific \cite{Perrone1975} \cite{kadalbajoo2010}.

Finite differences on nonuniform grids are routinely constructed by (i) choosing a monitor/weight g to control local resolution, (ii) equidistributing computational coordinates to obtain the grid, and (iii) applying finite-difference stencils on uneven spacing. This note collects those steps into a compact, reproducible workflow. We shall not propose a new mapping or new stencils. Rather, we gather standard formulas in one place, list practical conditions on g, and provide an application example to the 2D Schrödinger equation in a Harmonic Oscillator, comparing a uniform and a weight-driven mesh at fixed matrix size.

We review the 1-D construction $\xi=S/S_{total}$ with $S(x)=\int_a^x 1/g(s)\,ds$ and its separable multi-dimensional analogue, including the frequently used harmonic/variable-diffusion mapping $\nabla_{\xi}\!\cdot((1/g)\nabla_{\xi}x)=0$. For derivatives, we use the three-point central approximations on uneven spacing (equal to the standard central differences on a uniform grid), which we write explicitly for convenience.

\subsection{Relevant Works}

This note builds on the well-established existing work on uneven-grid finite-differences, equidistribution and monitor functions and harmonic grid generation. Closed-form and algorithmic generation of finite-difference weights on arbitrarily space nodes is discussed e.g. in \cite{Fornberg1988} \cite{Fornberg1998}. Adaptive ad moving-mesh methods are discussed at length e.g. in Huang and Russell's monograph in Reference \cite{Huang2011}. As for harmonic grid generation, the multidimensional mapping $\nabla_{\xi}\!\cdot((1/g)\nabla_{\xi}x)=0$ which will be used in this work is a special case of a Winslow-type mapping with scalar diffusion. Overviews include Thompson, Warsi and Mastin's book \cite{Thompson1985} and later exposition connecting harmonic maps and grid quality \cite{Farrashkhalvat2003}.

\section{One-Dimensional Case}
\subsection{Mesh Generation}
Before deriving the finite differences in a variable mesh, it is useful to define how this is constructed.

Firstly, we define a weighting function $g(x)$. This function gives the relative step size as a function of the position in parameter space. We say relative because due to the way we will construct the mesh, this function will always be normalized to 1, so, for instance, $g(x)=k$, where k is a constant, will yield the exact same uniformly spaced grid for any k. This function simply describes the relative density of the mesh. The function's co-domain should lie in $(0,1]$ or any linear scaling of that. Positions in the parameter space where the value of $g(x)$ is lower (higher) will correspond to areas with a finer (coarser) mesh. We will see later that this function essentially generates the Cumulative Distribution Function (CDF) of the mesh spacing.

Let us now say that we are trying to generate a grid of points $x_i$ over the interval $[a,b]$, using the previously defined weighting function $g(x)$, and dividing the interval into N segments. We shall employ a mapping that transforms a uniformly spaced computational domain into our variably spaced physical domain. We hence define the computational coordinate $\xi\in[0,1]$ over the interval $[a,b]$:
\begin{equation}
    \xi = \frac{S(x)}{S_{\text{total}}},
    \label{xi_1D}
\end{equation}
where S(x) is the CDF of the spacing over $[a,x]$ defined from the weighting function $g(x)$ as:
\begin{equation}
    S(x) = \int_{a}^{x} \frac{1}{g(s)} \, ds,
    \label{s_function_1d}
\end{equation}
and $S_{\text{total}}$ is the total "scaled" length of the interval:
\begin{equation}
    S_{\text{total}} = S(b) = \int_{a}^{b} \frac{1}{g(s)} \, ds.
\end{equation}
The reasoning for this choice will become clearer later upon extrapolation to higher dimensions. For now, we note that it establishes a diffeomorphism between the computational and physical spaces, and we can freely transform between the two via this function. We can now generate uniform grid points in $\xi$-space by dividing its domain into N equal segments as:
\begin{equation}
    \xi_i = \frac{i}{N}, \quad i = 0, 1, 2, \dots, N,
\end{equation}
and map these points from the computational grid onto the physical grid via the relation:
\begin{equation}
    S(x_i)=\xi_i\cdot S_{\text{total}}\Rightarrow x_i=S^{-1}(\xi_i\cdot S_{\text{total}}),
    \label{s_1d}
\end{equation}
where $S(x)$ might not always have an analytical inverse, but it can be easily found via numerical methods in a practical implementations. 
 
\subsection{Finite Differences}
Let $x_{i-1},x_i,x_{i+1}$ be three consecutive points, let $h_{i-1}=x_i-x_{i-1}$ be the distance between $x_{i-1}$ and $x_i$ and let $h_i=x_{i+1}-x_i$ be the distance between $x_i$ and $x_{i+1}$. We want to derive the finite differences that approximate the derivatives of a function $f(x)$ in a discrete mesh using a central approximation at a point $x_i$. Here, we derive the first- and second-order derivatives, but the same principle applies further. We want to find the coefficients $A_i,B_i,C_i$ where i=1,2, such that:
\begin{equation}
    f'_i \approx A_1 f_{i-1} + B_1 f_i + C_1 f_{i+1}
\end{equation}
\begin{equation}
    f''_i \approx A_2 f_{i-1} + B_2 f_i + C_2 f_{i+1}
\end{equation}
We then substitute $f_{i-1}$ and $f_{i+1}$ with their respective Taylor series expansions around $x_i$:
\begin{equation}
    f_{i-1} = f_i - h_{i-1} f'_i + \frac{h_{i-1}^2}{2} f''_i - \frac{h_{i-1}^3}{6} f'''_i + \cdots
\end{equation}
\begin{equation}
    f_{i+1} = f_i + h_i f'_i + \frac{h_i^2}{2} f''_i + \frac{h_i^3}{6} f'''_i + \cdots
\end{equation}
Solving, we find that the coefficients are
\begin{equation}
\begin{aligned}
A_1 & = - \frac{h_i}{h_{i-1} (h_i + h_{i-1})} \\
B_1 & = \frac{h_i^2 - h_{i-1}^2}{h_i h_{i-1} (h_i + h_{i-1})} \\
C_1 & = \frac{h_{i-1}}{h_i (h_i + h_{i-1})}
\end{aligned}
\end{equation}
and
\begin{equation}
\begin{aligned}
A_2 & = \frac{2}{h_{i-1} (h_{i-1} + h_i)} \\
B_2 & = - \frac{2}{h_{i-1} h_i} \\
C_2 & = \frac{2}{h_i (h_{i-1} + h_i)}
\end{aligned}
\end{equation}
Therefore, the general form of the finite differences with variable mesh becomes, for the first derivative:
\begin{equation}
\begin{split}
    f'_i \approx - \frac{h_i}{h_{i-1} (h_i + h_{i-1})} & f_{i-1} + \frac{h_i^2 - h_{i-1}^2}{h_i h_{i-1} (h_i + h_{i-1})} f_i + \\
    & +\frac{h_{i-1}}{h_i (h_i + h_{i-1})} f_{i+1}
\end{split}
\end{equation}
and for the second derivative:
\begin{equation}
\begin{split}
    f''_i \approx \frac{2}{h_{i-1} (h_{i-1} + h_i)} f_{i-1} - & \frac{2}{h_{i-1} h_i} f_i + \\
    & +\frac{2}{h_i (h_{i-1} + h_i)} f_{i+1}
\end{split}
\end{equation}
One can easily see that in the case of a uniform grid, $h_{i-1}=h_i\equiv h$, and hence
\begin{equation}
    f'_i \approx - \frac{1}{2h} f_{i-1} + 0 \cdot f_i + \frac{1}{2h} f_{i+1} = \frac{f_{i+1} - f_{i-1}}{2h}
\end{equation}
\begin{equation}
    f''_i \approx \frac{1}{h^2} f_{i-1} - \frac{2}{h^2} f_i + \frac{1}{h^2} f_{i+1} = \frac{f_{i+1} - 2 f_i + f_{i-1}}{h^2}
\end{equation}
and we retrieve the known form for a constant-step mesh as in Equations \ref{finitediff_order1} and \ref{finitediff_order2}.

\section{Generalization to Higher Dimensions}
\subsection{Mesh Generation}
\subsubsection{General weight function}

In the case of multiple dimensions, analogously to the method previously described for a single dimension, our aim is to map a uniformly spaced computational domain $\mathbf{\xi}=(\xi_1,\xi_2,...,\xi_n) \in [0, 1]^n$ to the physical domain $\mathbf{x}=(x_1,x_2,...,x_n)$, according to the weight function $g(\mathbf{x})=g(x_1,x_2,...,x_n)$. We can define the map between the computational and physical spaces by considering the divergence of a vector field in computational space:
\begin{equation}
    \nabla_{\boldsymbol{\xi}} \cdot \left( \mathbf{P} \nabla_{\boldsymbol{\xi}} \mathbf{x} \right) = \mathbf{0},
\end{equation}
where $\mathbf{P}$ is a tensor that incorporates the grid density, and in our case it is generated by the weight function as
\begin{equation}
    \mathbf{P} = \frac{1}{g(\mathbf{x})} \mathbf{I},
\end{equation}
where $\mathbf{I}$ is the identity tensor. This yields the form of the grid generation equations given a fully general $n$-dimensional $g(\mathbf{x})$
\begin{equation}
    \nabla_{\boldsymbol{\xi}} \cdot \left( \frac{1}{g(\mathbf{x})} \nabla_{\boldsymbol{\xi}} \mathbf{x} \right) = 0,
\end{equation}
which is a system of $n$ coupled PDEs for $\mathbf{x}(\mathbf{\xi})$. Note that for $n=1$, this reduces to a single equation
\begin{equation}
    \frac{d}{d\xi} \left( \frac{1}{g(x)} \frac{dx}{d\xi} \right) = 0
\end{equation}
\begin{equation}
    \frac{dx}{d\xi} = C \cdot g(x) \Rightarrow \frac{d\xi}{dx} = \frac{C}{g(x)},
\end{equation}
where by integrating and identifying C as $1/S_{\text{total}}$, we can easily retrieve the form in Equation \ref{xi_1D}. To ensure that the diffeomorphism between computational and physical spaces holds, we must set that $g(\mathbf{x})$ should be integrable everywhere in the domain, different from zero for any $\mathbf{x}$ and be bounded, such that $1/g(\mathbf{x})\neq 0$.

\subsubsection{Separable weight function}

In the case where $g(\textbf{x})$ is separable:
\begin{equation}
    g(\mathbf{x}) = \prod_{i=1}^n g_i(x_i),
\end{equation}
the grid generation equations can be decoupled into a set of $n$ independent PDEs
\begin{equation}
    \frac{\partial}{\partial \xi_i} \left( \frac{1}{g_i(x_i)} \frac{\partial x_i}{\partial \xi_i} \right ) = 0.
\end{equation}
We can simplify this further as
\begin{equation}
    \frac{1}{g_i(x_i)} \frac{\partial x_i}{\partial \xi_i} = \text{k}_i,
\end{equation}
where $k_i$ are constants. Finally, integrating, we arrive at the relation
\begin{equation}
    x_i = S_i^{-1} (\xi_i \cdot k_i)
\end{equation}
where $S_i(x_i)$ is completely analogous to the definition introduced in Equation \ref{s_function_1d}
\begin{equation}
    S_i(x_i) = \int_{a_i}^{x_i} \frac{1}{g_i(s)} \, ds.
\end{equation}
We note that separability allows for an independent and one-dimensional treatment of each dimension, which greatly simplifies the mesh generation problem. For that reason, it might be desirable in most cases to use weight functions of this type to speed up and simplify the grid generation. It is easy to see that $k_i$ can be identified with the normalization constant of the weight function in each dimension $i$, and hence we fully retrieve our treatment for the one-dimensional case as described in Equation \ref{s_1d}.

\subsection{Finite Differences}

For ease of notation, let's define
\begin{equation}
    \mathbf{x}=(x_1,x_2,x_3,...)\equiv(x,y,z,...).
\end{equation}
Consider that we have generated the grid points $\mathbf{x}_{ijk...}=(x_i, y_j,z_k,..)$, where $i,j,k,...$ index the grid in the $x,y,z,...$ directions. The step sizes in the direction $x$ are, analogously to the one-dimensional case:
\begin{equation}
    h_{i-1}^x = x_i - x_{i-1}
\end{equation}

\begin{equation}
    h_i^x = x_{i+1} - x_i
\end{equation}
And equivalently for all other dimensions. Following a process analogous to the one shown for the one-dimensional case, we find that the finite difference approximations to the first and second derivatives are given by Equations \ref{1st_order_higherdim} and \ref{2nd_order_higherdim} respectively. Note that these Equations give the form for the partial derivatives with respect to $x$, but this is perfectly extendable to any other variable by simply replacing $x$ with the variable of interest and affecting the respective index instead of $i$.
\begin{widetext}
\begin{equation}
    \left( \frac{\partial f}{\partial x} \right)_{i, j, k, \dots} \approx -\frac{h_i^x}{h_{i-1}^x(h_i^x + h_{i-1}^x)} f_{i-1, j, k,\dots} + \frac{(h_i^x)^2 - (h_{i-1}^x)^2}{h_i^x h_{i-1}^x (h_i^x + h_{i-1}^x)} f_{i, j, k,\dots} + \frac{h_{i-1}^x}{h_i^x(h_i^x + h_{i-1}^x)} f_{i+1, j, k,\dots}
    \label{1st_order_higherdim}
\end{equation}
\begin{equation}
    \left( \frac{\partial^2 f}{\partial x^2} \right)_{i, j, k,\dots} \approx \frac{2}{h_{i-1}^x(h_i^x + h_{i-1}^x)} f_{i-1, j, k,\dots} - \frac{2}{h_i^x h_{i-1}^x} f_{i, j, k, \dots} + \frac{2}{h_i^x(h_i^x + h_{i-1}^x)} f_{i+1, j, k,\dots}
    \label{2nd_order_higherdim}
\end{equation}
\end{widetext}

\section{Example of Implementation}

\begin{figure*}[t!]
    \centering
    \begin{subfigure}[t]{0.49\textwidth}
        \centering
        \includegraphics[width=0.8\textwidth]{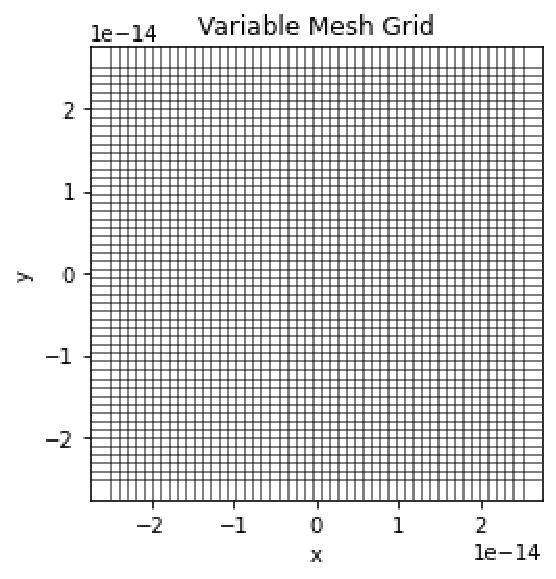}
        \caption{Uniform mesh}
    \end{subfigure}%
    ~ 
    \begin{subfigure}[t]{0.49\textwidth}
        \centering
        \includegraphics[width=0.8\textwidth]{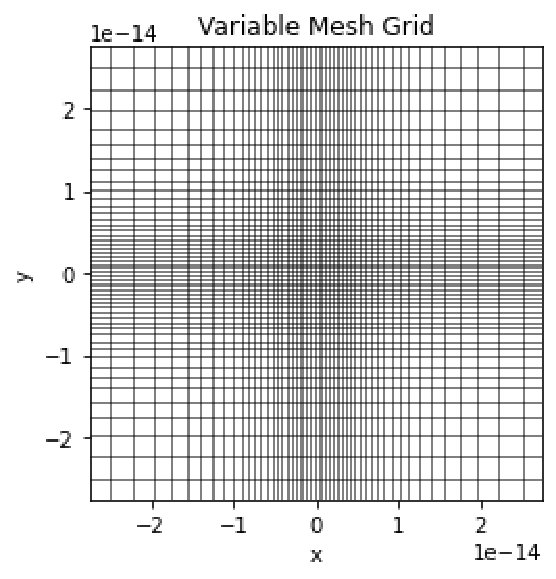}
        \caption{Variable mesh}
        \label{Variable_mesh}
    \end{subfigure}
    \caption{Mesh construction comparison. Both $x$ and $y$ are in meters.}
    \label{Mesh_comparison}
\end{figure*}

We shall use a grid of $N_x\times N_y=50\times50$ points in both cases, meaning that the Hamiltonian matrix of this eigenvalue problem has the same size in both cases, so the problems are computationally equivalent. The weight functions used were:
\begin{equation}
    g_x(x)=1-0.9\exp\left[-\left(\frac{x}{b_x-a_x}\right)^2\right],
\end{equation}
\begin{equation}
    g_y(y)=1-0.9\exp\left[-\left(\frac{y}{b_y-a_y}\right)^2\right].
\end{equation}
We hold $N_x\times N_y$ fixed so the Hamiltonian size is identical across meshes. Therefore, improvements in the resolution of the eigenvalues come purely from point redistribution according to the weight functions above. The resulting mesh is shown in Figure \ref{Variable_mesh}, next to the uniform mesh, for comparison. 

The two-dimensional time-independent Schrödinger equation is given by \cite{Schrdinger1926}:
\begin{equation}
    E\psi(\mathbf{x})=-\frac{\hbar^2}{2m}\left(\frac{\partial^2}{\partial x^2}+\frac{\partial^2}{\partial y^2}\right)\psi(\mathbf{x})+V(\mathbf{x})\psi(\mathbf{x})
    \label{schrodinger2d}
\end{equation}
Which can be formulated in the form of an eigenvalue problem $E\psi=\hat{H}\psi$, by defining a Hamiltonian operator
\begin{equation}
    \hat{H}=-\frac{\hbar^2}{2m}\left(\frac{\partial^2}{\partial x^2}+\frac{\partial^2}{\partial y^2}\right)+V(\mathbf{x}).
\end{equation}
For a harmonic oscillator, we use the potential
\begin{equation}
    V(x,y)=\frac{1}{2}m\omega(x^2+y^2).
    \label{potential_ho}
\end{equation}
For this example, we shall discretize the two-dimensional phase space $x,y\in[-25,25]\text{ fm}$ in 50 steps. As for the potential, we use the form described in Equation \ref{potential_ho} using $m=m_p$ and $\omega\hbar=10\text{ MeV}=1.6022\cdot 10^{-12}\text{ J}$. 
We can discretize Equation \ref{schrodinger2d} using the finite difference scheme described in Equation \ref{2nd_order_higherdim}.
\begin{widetext}
\begin{equation}
\begin{aligned}
- \frac{\hbar^2}{2m} &\left( \frac{2}{h_{i-1}^x (h_{i-1}^x + h_i^x)} \psi_{i-1, j} - \frac{2}{h_{i-1}^x h_i^x} \psi_{i, j} + \frac{2}{h_i^x (h_{i-1}^x + h_i^x)} \psi_{i+1, j} \right. \\
&+ \left. \frac{2}{h_{j-1}^y (h_{j-1}^y + h_j^y)} \psi_{i, j-1} - \frac{2}{h_{j-1}^y h_j^y} \psi_{i, j} + \frac{2}{h_j^y (h_{j-1}^y + h_j^y)} \psi_{i, j+1} \right ) + V_{i,j} \psi_{i,j} = E \psi_{i,j}.
\label{finite_differences_sch2d}
\end{aligned}        
\end{equation}
\end{widetext}
By treating $\psi_{ij}$ as a $(i\times j)\times 1$ matrix, we separate the left-hand-side into an $(i\times j)\times(i \times j)$ matrix multiplying by $\psi_{ij}$, and arrive at the well-known eigenvalue problem form
\begin{equation}
    \hat{H}_{ij}\psi_{ij}=E\psi_{ij},
\end{equation}
which was then solved using Python.

\begin{figure*}[htbp]
\centering
\begin{tabular}{c@{\hspace{1.5cm}}c}
    % First Row
    \includegraphics[width=\figsize\textwidth]{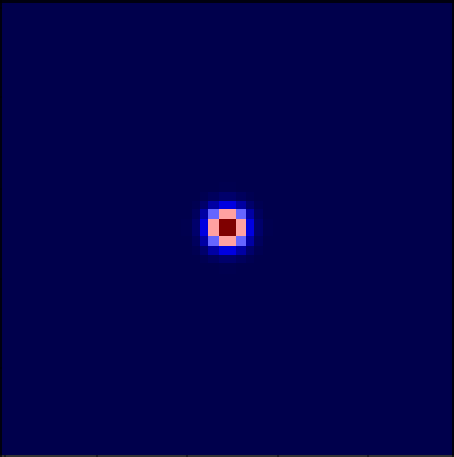} &
    \includegraphics[width=\figsize\textwidth]{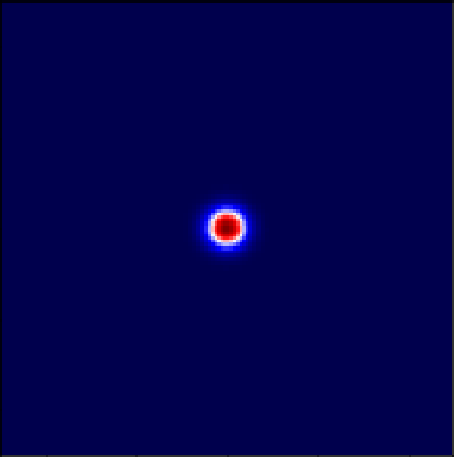} \\
    \multicolumn{2}{c}{(1)} \\[0.5ex]
    
    % Second Row
    \includegraphics[width=\figsize\textwidth]{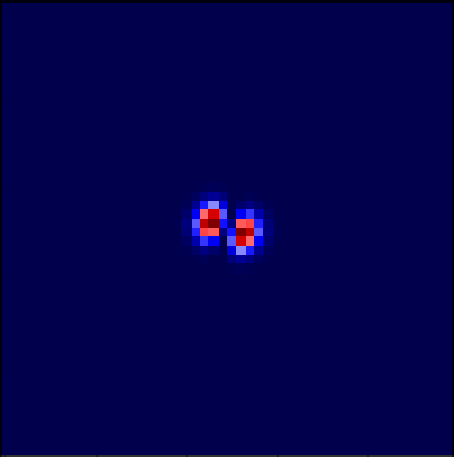} &
    \includegraphics[width=\figsize\textwidth]{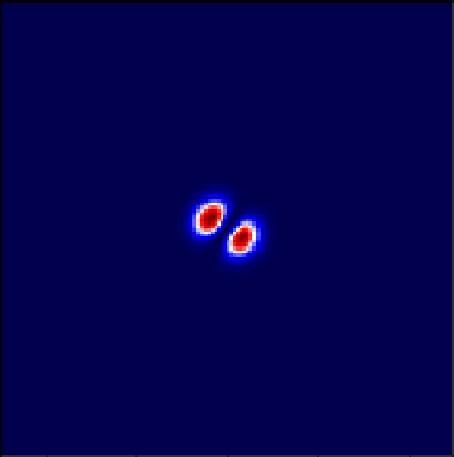} \\
    \multicolumn{2}{c}{(2)} \\[0.5ex]
    
    % Third Row
    \includegraphics[width=\figsize\textwidth]{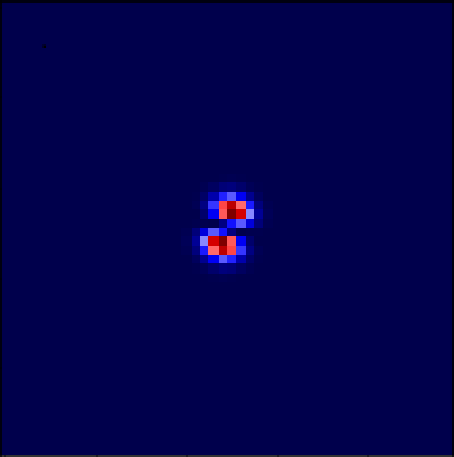} &
    \includegraphics[width=\figsize\textwidth]{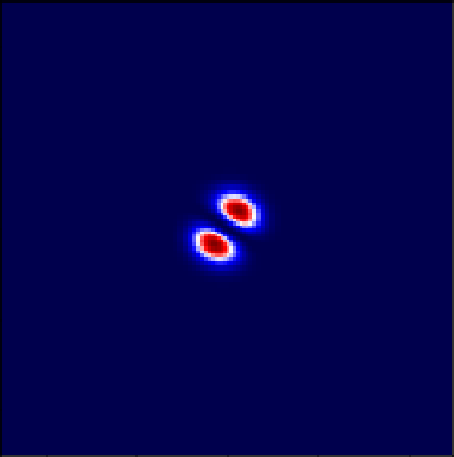} \\
    \multicolumn{2}{c}{(3)} \\[0.5ex]
    
    % Fourth Row
    \includegraphics[width=\figsize\textwidth]{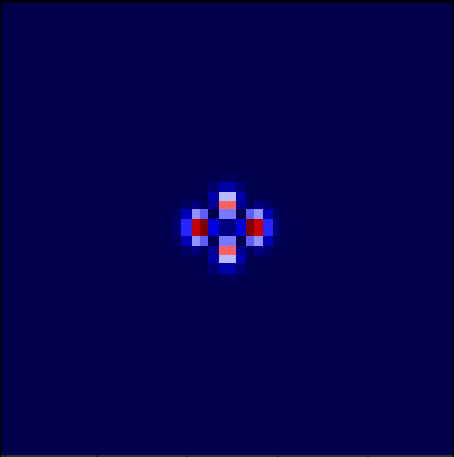} &
    \includegraphics[width=\figsize\textwidth]{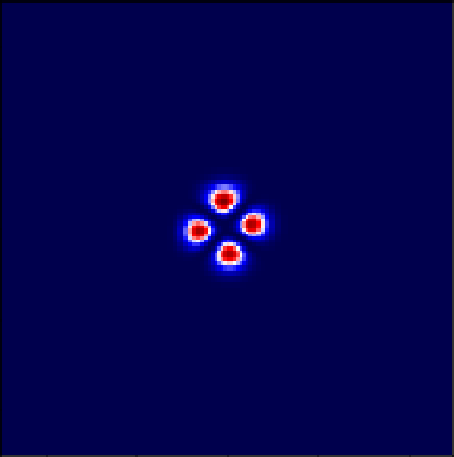} \\
    \multicolumn{2}{c}{(4)} \\[0.5ex]
    
    % Fifth Row
    \includegraphics[width=\figsize\textwidth]{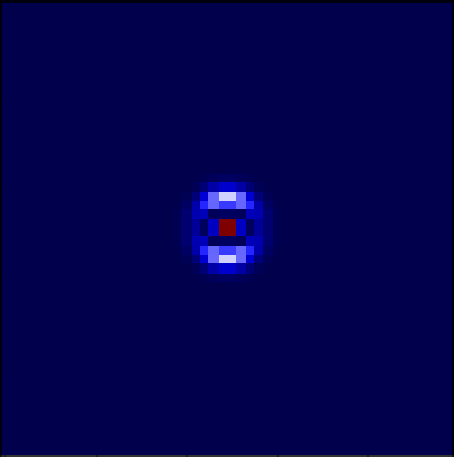} &
    \includegraphics[width=\figsize\textwidth]{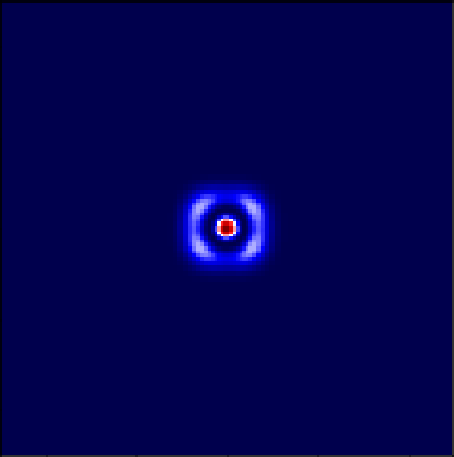} \\
    \multicolumn{2}{c}{(5)} \\
\end{tabular}
\caption{First five lowest-energy eigenfunctions of the 2D Harmonic Oscillator. Both meshes have the same average step size which are rearranged through the action of weight functions. Left: uniform mesh; Right: variable mesh;}
\label{fig:my_figure}
\end{figure*}

\section{Conclusions}

We have derived a formalism that allows the application of the Finite Difference Method using a variable mesh through a relative mesh spacing weight function $g(\mathbf{x})$. This formalism is presented for both the one-dimensional and arbitrary-dimensional case. We have shown an example of application by solving the two-dimensional Schrödinger equation using both uniform and variable mesh, with the same average step size, observing how this formalism can be a useful tool to better resolve the eigenfunctions of this Hamiltonian, which are very localized at the center of the parameter space, without necessarily increasing the computational cost, since ultimately the total size of the Hamiltonian matrix whose eigenvalues and eigenfunctions we are solving for is the same for both. This note consolidates existing techniques into a minimal recipe. We do not claim new stencils or a new grid-generation PDE.

\bibliography{apssamp}% Produces the bibliography via BibTeX.

\begin{acknowledgments}
The author would like to thank Chong Qi and Rodrigo Arouca for the fruitful discussions on applications.
\end{acknowledgments}

\end{document}